\documentclass[reqno,oneside,english]{amsart}
\usepackage[T1]{fontenc}
\usepackage[latin9]{inputenc}
\usepackage{amsmath}
\usepackage{amsthm}
\usepackage[pdftex]{graphicx}
\usepackage{geometry}
\usepackage{amssymb}

\makeatletter
\numberwithin{equation}{section}
\numberwithin{figure}{section}

\makeatother

\usepackage{babel}
\begin{document}

\title[\tiny{The Laplacian on the unit square in a self-similar manner}]{The Laplacian on the unit square in a self-similar manner}
\author{Hua Qiu}
\address{Department of Mathematics, Nanjing University, Nanjing, 210093, P. R. China.}
\curraddr{} \email{huaqiu@nju.edu.cn}
\thanks{}

\author{Haoran Tian$^*$}
\address{School of Physics, Nanjing University, Nanjing, 210093, P.R. China.}
\curraddr{} \email{hrtian@hotmail.com}
\thanks{$^*$ Corresponding author. }
\thanks{The research of the first author was supported by the Nature Science Foundation of China, Grant 11471157.}



\date{}

\dedicatory{}
\begin{abstract}
In this paper, we show how to construct the standard Laplacian on the unit square  in a self-similar manner. We rewrite the familiar mean value property of planar harmonic functions in terms of averge values on small squares, from which we could know how the planar self-similar resistance form and the Laplacian look like. This approach combines the constructive limit-of-difference-quotients method of Kigami for p.c.f. self-similar sets and the method of averages introduced by Kusuoka and Zhou for the Sierpinski carpet.
\end{abstract}
\maketitle

\section{Introduction}

Let $K$ denote a \textit{self-similar set} which is non-empty compact and generated by a finite family of  contraction
 similarity mappings  $\{F_i\}_{i=1,\dots,N}$ on $\mathbb{R}^n$ such that

\begin{equation*}
K=\bigcup_{i=1}^NF_iK.
\end{equation*}
If $w=(w_1,\cdots,w_m)$ is a finite word with each $w_j\in\{1,\cdots,N\}$, we define the mapping $$F_w=F_{w_1}\circ\cdots\circ F_{w_m},$$
and call  $F_wK$ a $m$-\textit{cell}  of $K$.

$K$ is called \textit{post-critically finite (p.c.f.)} if $K$ is connected, and there exists a finite set $V_0\subseteq K$ called the \textit{boundary}, such that
\begin{equation*}
\begin{aligned}
F_wK\cap F_{w'}K&\subseteq F_wV_0\cap F_{w'}V_0, 
\quad(F_wV_0\cap F_{w'}V_0)\cap V_0=\emptyset,\\ 
&\text{for} \ w\ne w' \  \text{with} \ |w|=|w'|,
\end{aligned}
\end{equation*}
where $|w|$ is the length of $w$.
Moreover, we require that each boundary point is the fixed point of one of the mappings in $\{F_i\}_{i=1,\dots,N}$. Without loss of generality we  assume that $V_0=\{q_1,\dots,q_{N_0}\}$ for $N_0\leq N$ with
\begin{equation*}
F_iq_i=q_i, \text{ for } i=1,\dots,N_0.
\end{equation*}

A more general definition of p.c.f. self-similar sets can be found in [2,3] introduced by Kigami. The unit interval(\textit{I}) and the Sierpinski gasket(\textit{SG}) are two typical examples. However, the unit square ($S$, which we mainly discuss in this paper) and the Sierpinski carpet(\textit{SC}) are non-p.c.f. self-similar sets. 

An analytic construction of a Laplacian on \textit{SG} was given by Kigami [1] as a \textit{renormalized limit} of difference quotients, which he later extended to p.c.f. self-similar sets [2]. Please also refer to [6] to find a detailed introduction to this topic.  At about the same time, Kusuoka and Zhou [4] developed a \textit{method of average} for defining a Laplacian on $SC$, which uses average values of functions over cells rather than pointwise values in the definition. This method was later proved  to be equivalent to Kigami's definition for \textit{symmetric Laplacian}(with respect to the standard self-similar measure) by Strichartz [5] on \textit{SG}, which brings the possibility to define Laplacians on more general self-similar sets. 

Recall that $SG$ is the invariant set of a  family of $3$ contraction mappings $F_ix=\frac{1}{3}(x+q_i), i=1,2,3$, where $q_1,q_2,q_3$ are the vertices of an equilateral triangle in $\mathbb{R}^2$.  Let $\mu$ be the standard \textit{regular probability measure} with weights $(\frac{1}{3},\frac{1}{3},\frac{1}{3})$ on \textit{SG}. Define the average value for a function $f$ on the cell $F_wSG$ as 
\begin{equation*}
B_w(f)=\frac{1}{\mu(F_wSG)}\int_{F_wSG}fd\mu.
\end{equation*}
It is easy to check that 
$$
B_w(f)=\frac{1}{3}\sum_{i=1}^3B_{wi}(f),
$$
where $wi$ denotes $(w_1,\cdots,w_m,i)$.

According to [5], any harmonic function $h$ on \textit{SG} satisfies the \textit{mean value property}
\begin{equation*}
B_w(h)=\frac{1}{3}\sum_{w'\sim_mw}B_{w'}(h)
\end{equation*}
for any finite word $w$ with $F_wSG\cap V_0=\emptyset$, where $|w'|=|w|$ and $w'\sim_m w$ means that $F_{w'}SG\cap F_wSG$ is non-empty. 

The resistance form $\mathcal{E}(\cdot,\cdot)$ on \textit{SG} is defined as
\begin{equation*}
\mathcal{E}(f,g)=\lim_{m\to\infty}\mathcal{E}_m(f,g),
\end{equation*}
with
\begin{equation*}
\mathcal{E}_m(f,g)=\frac{3}{2}r^{-m}\sum_{w'\sim_mw}\big(B_{w'}(f)-B_w(f)\big)\big(B_{w'}(g)-B_w(g)\big),
\end{equation*}
where $r=\frac{3}{5}$ is the renormalization factor.

The symmetric Laplacian was proved to be the limit of a sequence of discrete Laplacians $\Delta_m$, in the sense that 
\begin{equation*}
\Delta f(x)=\lim_{m\to\infty}\Delta_mf(x)
\end{equation*}
uniformly for every $f$ in the domain of $\Delta$, where $\Delta_m$ is defined as
\begin{equation*}
\Delta_mf(x)=\frac{9}{2}r^{-m}3^m\bigg(\frac{1}{3}\sum_{w'\sim_mw}B_{w'}(f)-B_w(f)\bigg), \text{ for } x\in F_wK,|w|=m,
\end{equation*}
for each integer $m\geq 0$.

In essential, the method of average studies first the analysis on \textit{cell graphs} $\Gamma_m$ whose vertices are the words $w$ of length $m$, and edge relation, denoted by $w\sim_m w'$ is defined by the condition that $F_wSG\cap F_{w'}SG\neq\emptyset$, then passes the approximation to the limit as $m\rightarrow\infty$ to establish the analysis on $SG$.

In this paper, we will investigate that to what extent can we extend Kusuoka and Zhou's method to the planar unit square $S$, which could be viewed as a non-p.c.f. self-similar set, generated by a family of contraction mappings $\{F_i\}_{i=1}^4$ on $\mathbb{R}^2$ with
\begin{equation*}
F_i(x)=\frac{1}{2}(x+q_i), \ i=1,2,3,4,
\end{equation*}
and
\begin{equation*}
q_1=(0,0), \ q_2=(1,0), \ q_3=(1,1), \ q_4=(0,1).
\end{equation*}
In particular, the mean value property for planar harmonic functions on cell graphs(several choices) will not hold exactly. However, it will be close to holding exactly for large $m$. We will provide a sharp estimate for it.  The analysis established in this paper provides a self-similar viewpoint of the familiar classical analysis on $\mathbb{R}^2$. It will be nice to be able to transform the classical analytical results via  this self-similar construction, which maybe will bring insights to the analysis on $SC$. 

\section{Mean value property of harmonic functions}

From now on, we use $\mu$ to denote the \textit{Lebesgue measure} on $\mathbb{R}^2$. For $x=(x_1,x_2)\in\mathbb{R}^2$, $l\in\mathbb{R}^+$, denote
\begin{equation}
\mathcal{D}(x,l)=\{(\xi_1,\xi_2)\in\mathbb{R}^2|x_1-\frac{l}{2}\leq \xi_1\leq x_1+\frac{l}{2},x_2-\frac{l}{2}\leq \xi_2\leq x_2+\frac{l}{2}\}
\end{equation}
 a $l$-\textit{square} centered at $x$. For a function $f$ integrable on $\mathcal{D}(x,l)$, write
\begin{equation}
I(f,x,l)=\frac{1}{l^2}\int_{\mathcal{D}(x,l)}fd\mu.
\end{equation}

\noindent\textbf{Definition 2.1.} \textit{Let $p=(p_1,p_2)$ be an integer pair with $0\leq p_1\leq p_2$. For two squares $\mathcal{D}(x,l)$ and $\mathcal{D}(x',l)$, write $\mathcal{D}(x,l)\sim_p \mathcal{D}(x',l)$ if
\begin{equation}
\begin{aligned}
&|x_1-x_1'|=p_1l, \ |x_2-x_2'|=p_2l,\ \text{or}\\
&|x_1-x_1'|=p_2l, \ |x_2-x_2'|=p_1l,
\end{aligned}
\end{equation}
call them $p$-neighbors. 
}

Suppose $f$ is integrable on every $l$-squares $p$-neighbor to $\mathcal{D}(x,l)$, denote
\begin{equation}
I_p(f,x,l)=\sum_{\tiny\substack{
\mathcal{D}(x',l) \\ 
\sim_p \mathcal{D}(x,l)}
}\!\!\!\!I(f,x',l).
\end{equation}

For simplicity, denote a constant $c_p$ as
\begin{equation}
c_p \ = \
\begin{cases}
\begin{aligned}
\frac{1}{8}&,\text{ if } p_1=p_2=0, \\
\frac{1}{2}&,\text{ if } p_1=0,p_2\ne 0, \text{ or } p_1=p_2\ne 0, \\
1&,\text{ if }0<p_1<p_2.
\end{aligned}
\end{cases}
\end{equation}
It is easy to check that $8c_p$ is the number of $p$-neighbors of any fixed $\mathcal{D}(x,l)$.

\noindent\textbf{Lemma 2.2.} \textit{Let $x\in\mathbb{R}^2$, $n$ be a non-negative integer. For functions $f_n^x$, $g_n^x$ defined by
\begin{equation}
\begin{aligned}
&f_n^x(\xi)=\sum_{j=0}^{[\frac{n}{2}]}(-1)^j\frac{(\xi_1-x_{1})^{n-2j}(\xi_2-x_{2})^{2j}}{(n-2j)!(2j)!},\\
&g_n^x(\xi)=\sum_{j=0}^{[\frac{n-1}{2}]}(-1)^j\frac{(\xi_1-x_{1})^{n-2j-1}(\xi_2-x_{2})^{2j+1}}{(n-2j-1)!(2j+1)!},
\end{aligned}
\end{equation}
and integer pair $p=(p_1,p_2)$ with $0\leq p_1\leq p_2$, we have
\begin{equation}
\begin{aligned}
&I_p(f_n^x,x,l)=
\frac{1}{(4k+2)!}T_p^{(k)}(\frac{l}{2})^{4k}, \text{ if } \ 4\mid n \text{ with } n=4k,\\
&I_p(f_n^x,x,l)=0, \text{ if } \ 4\nmid n,\\
&I_p(g_n^x,x,l)=0,
\end{aligned}
\end{equation}
where constants $T_p^{(k)}$ are
\begin{equation}
\begin{aligned}
T_p^{(k)}= \  &2c_p\Bigg(\left((2p_1+1)^2+(2p_2+1)^2\right)^{2k+1}\sin\left((4k+2)\arctan\frac{2p_2+1}{2p_1+1}\right) \\
&\quad-\left((2p_1-1)^2+(2p_2+1)^2\right)^{2k+1}\sin\left((4k+2)\arctan\frac{2p_2+1}{2p_1-1}\right) \\
&\quad+\left((2p_1-1)^2+(2p_2-1)^2\right)^{2k+1}\sin\left((4k+2)\arctan\frac{2p_2-1}{2p_1-1}\right) \\ &\quad-\left((2p_1+1)^2+(2p_2-1)^2\right)^{2k+1}\sin\left((4k+2)\arctan\frac{2p_2-1}{2p_1+1}\right) \Bigg),
\end{aligned}
\end{equation}
with the estimates that
\begin{equation}
|T_p^{(k)}|\leq 8c_p\|2p+(1,1)\|^{4k+2}.
\end{equation}
}

\noindent\textit{Proof.} By direct calculation, we obtain
\begin{equation*}
\begin{aligned}
I_p(f_n^x,x,l)=\\
\frac{c_p}{l^2}\sum_{j=0}^{[\frac{n}{2}]}&\frac{(-1)^j}{(n-2j+1)!(2j+1)!}\Bigg(\Big(\big((2p_1+1)\frac{l}{2}\big)^{n-2j+1}+\big(-(2p_1-1)\frac{l}{2}\big)^{n-2j+1}\\
-\big(&(2p_1-1)\frac{l}{2}\big)^{n-2j+1}-\big(-(2p_1+1)\frac{l}{2}\big)^{n-2j+1}\Big)\Big(\big((2p_2+1)\frac{l}{2}\big)^{2j+1}\\
+\big(&-(2p_2-1)\frac{l}{2}\big)^{2j+1}-\big((2p_2-1)\frac{l}{2}\big)^{2j+1}-\big(-(2p_2+1)\frac{l}{2}\big)^{2j+1}\Big)\\
+\Big(\big(&(2p_2+1)\frac{l}{2}\big)^{n-2j+1}+\big(-(2p_2-1)\frac{l}{2}\big)^{n-2j+1}-\big((2p_2-1)\frac{l}{2}\big)^{n-2j+1}\\
-\big(&-(2p_2+1)\frac{l}{2}\big)^{n-2j+1}\Big)\Big(\big((2p_1+1)\frac{l}{2}\big)^{2j+1}+\big(-(2p_1-1)\frac{l}{2}\big)^{2j+1}\\
-\big(&(2p_1-1)\frac{l}{2}\big)^{2j+1}-\big(-(2p_1+1)\frac{l}{2}\big)^{2j+1}\Big)\Bigg)
\end{aligned}
\end{equation*}
\begin{equation*}
\begin{aligned}
=c_p&\sum_{j=0}^{[\frac{n}{2}]}\frac{(-1)^j}{4(n-2j+1)!(2j+1)!}(\frac{l}{2})^n\Bigg(\Big(1-(-1)^{n-2j+1}\Big)\Big((1-(-1)^{2j+1}\Big)\\
&\Big((2p_1+1)^{n-2j+1}-(2p_1-1)^{n-2j+1}\Big)\Big((2p_2+1)^{2j+1}-(2p_2-1)^{2j+1}\Big)\\
+&\Big(1-(-1)^{n-2j+1}\Big)\Big((1-(-1)^{2j+1}\Big)\Big((2p_2+1)^{n-2j+1}-(2p_2-1)^{n-2j+1}\Big)\\
&\Big((2p_1+1)^{2j+1}-(2p_1-1)^{2j+1}\Big)\Bigg).
\end{aligned}
\end{equation*}
Obviously the summation is zero unless $n$ is even. If $n$ is even, let $n=2m$, then
\begin{equation*}
I_p(f_n^x,x,l)=\frac{1}{(2m+2)!}(\frac{l}{2})^{2m}c_p\sum_{j=0}^{m}(-1)^jC_{2m+2}^{2j+1}\gamma_p^{m-j,j},
\end{equation*}
with
\begin{equation*}
\begin{aligned}
\gamma_p^{m-j,j}=&\big((2p_1+1)^{2m-2j+1}-(2p_1-1)^{2m-2j+1}\big)\big((2p_2+1)^{2j+1}-(2p_2-1)^{2j+1}\big)\\
&+\big((2p_2+1)^{2m-2j+1}-(2p_2-1)^{2m-2j+1}\big)\big((2p_1+1)^{2j+1}-(2p_1-1)^{2j+1}\big).
\end{aligned}
\end{equation*}
Using the symmetry of $\sum_{j=0}^{m}(-1)^jC_{2m+2}^{2j+1}\gamma_p^{m-j,j}$, it is easy to check that the summation is zero unless $m$ is even. If m is even, let $m=2k$ ($n=2m=4k$), then we obtain that
\begin{equation*}
I_p(f_n^x,x,l)=\frac{1}{(4k+2)!}T_p^{(k)}(\frac{l}{2})^{4k},
\end{equation*}
and
\begin{equation*}
\begin{aligned}
T_p^{(k)}= \ &c_p\sum_{j=0}^{2k}(-1)^{j}C_{4k+2}^{2j+1}\gamma_p^{2k-j,j}\\
= \ &c_p\sum_{j=0}^{2k}(-1)^{j}C_{4k+2}^{2j+1}\Big(\big((2p_1+1)^{4k-2j+1}-(2p_1-1)^{4k-2j+1}\big)\big((2p_2+1)^{2j+1}-(2p_2-1)^{2j+1}\big)\\
&\qquad+\big((2p_2+1)^{4k-2j+1}-(2p_2-1)^{4k-2j+1}\big)\big((2p_1+1)^{2j+1}-(2p_1-1)^{2j+1}\big)\Big)\\
= \ &c_p\Big(Im\big(2p_1+1+i(2p_2+1)\big)^{4k+2}-Im\big(2p_1-1+i(2p_2+1)\big)^{4k+2}\\
& \ +Im\big(2p_1-1+i(2p_2-1)\big)^{4k+2}-Im\big(2p_1+1+i(2p_2-1)\big)^{4k+2}\\
& \ 
+Im\big(2p_2+1+i(2p_1+1)\big)^{4k+2}-Im\big(2p_2-1+i(2p_1+1)\big)^{4k+2}\\
& \ +Im\big(2p_2-1+i(2p_1-1)\big)^{4k+2}-Im\big(2p_2+1+i(2p_1-1)\big)^{4k+2}\Big).\\
\end{aligned}
\end{equation*}
Then an easy calculation yields (2.8) and (2.9).

Analogously, for $g_n^x$, a similar argument will gives
$
I_p(g_n^x,x,l)=0.
$
\hfill$\square$

In fact, $\{f_n^x\}$ and $\{g_n^x\}$ are \textit{polynomial harmonic functions} on $\mathbb{R}^2$, which form a "basis" of harmonic functions near $x$.

\noindent\textbf{Lemma 2.3.} \textit{Let $p$ be an integer pair as before, $\Omega$ be an open set in $\mathbb{R}^2$ and $x\in\Omega$. Then there exists a positive constant $l_{x,p}$, such that for any harmonic function $h$ on $\Omega$ and  $l<l_{x,p}$,
\begin{equation}
I_p(h,x,l)=8c_ph(x)+\sum_{k=1}^\infty\frac{1}{(4k+2)!}\frac{\partial^{4k}h}{\partial\xi_1^{4k}}(x)T_p^{(k)}(\frac{l}{2})^{4k},
\end{equation}
with the same constants $T_p^{(k)}$ as in (2.8).}

\noindent\textit{Proof.} Since $h$ is harmonic on $\Omega$, it is real analytic in $\Omega$. So we can expand $h$ near $x$ as
\begin{equation}
h(\xi)=\sum_{n=0}^\infty\sum_{i=0}^n\frac{\partial^nh}{\partial\xi_1^{n-i}\partial\xi_2^i}(x)\frac{(\xi_1-x_1)^{n-i}(\xi_2-x_2)^i}{(n-i)!i!}.
\end{equation}
Noticing that
\begin{equation*}
\Delta h(x)=\frac{\partial^2h}{\partial\xi_1^2}(x)+\frac{\partial^2h}{\partial\xi_2^2}(x)=0,
\end{equation*}
we have\begin{equation*}
\frac{\partial^nh}{\partial\xi_1^{n-i}\partial\xi_2^i}(x)=-\frac{\partial^nh}{\partial\xi_1^{n-i-2}\partial\xi_2^{i+2}}(x), \text{ for }  i=0,1,\dots,n\!-\!2.
\end{equation*}
By iterating we get
\begin{equation*}
\begin{aligned}
&\frac{\partial^nh}{\partial\xi_1^{n-2j}\partial\xi_2^{2j}}(x)=(-1)^j\frac{\partial^nh}{\partial\xi_1^n}(x), \text{ for } j=0,1,\dots,[\frac{n}{2}],\\
&\frac{\partial^nh}{\partial\xi_1^{n-2j-1}\partial\xi_2^{2j+1}}(x)=(-1)^j\frac{\partial^nh}{\partial\xi_1^{n-1}\partial\xi_2}(x), \text{ for } j=0,1,\dots,[\frac{n-1}{2}].
\end{aligned}
\end{equation*}
Then (2.11) gives
\begin{equation}
\begin{aligned}
h(\xi)
=&h(x)+\sum_{n=1}^\infty\sum_{j=0}^{[\frac{n}{2}]}\frac{\partial^nh}{\partial\xi_1^{n-2j}\partial\xi_2^{2j}}(x)\frac{(\xi_1-x_1)^{n-2j}(\xi_2-x_2)^{2j}}{(n-2j)!(2j)!}\\
&+\sum_{n=1}^\infty\sum_{j=0}^{[\frac{n-1}{2}]}\frac{\partial^nh}{\partial\xi_1^{n-2j-1}\partial\xi_2^{2j+1}}(x)\frac{(\xi_1-x_1)^{n-2j-1}(\xi_2-x_2)^{2j+1}}{(n-2j-1)!(2j+1)!}\\
=&h(x)+\sum_{n=1}^\infty\bigg(\frac{\partial^nh}{\partial\xi_1^n}(x)\sum_{j=0}^{[\frac{n}{2}]}(-1)^j\frac{(\xi_1-x_1)^{n-2j}(\xi_2-x_2)^{2j}}{(n-2j)!(2j)!}\bigg)\\
&+\sum_{n=1}^\infty\bigg(\frac{\partial^nh}{\partial\xi_1^{n-1}\partial\xi_2}(x)\sum_{j=0}^{[\frac{n-1}{2}]}(-1)^j\frac{(\xi_1-x_1)^{n-2j-1}(\xi_2-x_2)^{2j+1}}{(n-2j-1)!(2j+1)!}\bigg).\\
=&h(x)+\sum_{n=1}^\infty\bigg(\frac{\partial^nh}{\partial\xi_1^n}(x)f_n^x(\xi)+\frac{\partial^nh}{\partial\xi_1^{n-1}\partial\xi_2}(x)g_n^x(\xi)\bigg),
\end{aligned}
\end{equation}
where $f_n^x$ and $g_n^x$ are same as those in Lemma 2.2. 

Choose 
\begin{equation*}
l_{x,p}=\sup\{l>0| \text{ All } p\text{-neighbors of } \mathcal{D}(x,l) \text{ are contained in the convergence domain of (2.11)} \}.
\end{equation*}
Then for any $l<l_{x,p}$, we have 
\begin{equation}
I_p(h,x,l)=8c_ph(x)+\sum_{n=1}^\infty\bigg(\frac{\partial^nh}{\partial\xi_1^n}(x)I_p(f_n^x,x,l)+\frac{\partial^nh}{\partial\xi_1^{n-1}\partial\xi_2}(x)I_p(g_n^x,x,l)\bigg).
\end{equation}
By using Lemma 2.2, this gives (2.10).\hfill$\square$

The following is the \textit{mean value property} of planar harmonic functions in terms of average values on $l$-squares.

\noindent\textbf{Theorem 2.4.} \textit{Suppose $\mathcal{P}$ is a finite set consisting of integer pairs $p=(p_1,p_2)$ with $0\leq p_1\leq p_2$ and $p_2\ne 0$, and  $\{A_p\}_ {p\in\mathcal{P}}$ is a collection of real numbers satisfies
\begin{equation}
8\sum_{p\in\mathcal{P}}c_pA_p=1.
\end{equation}
Let $\Omega$ be an open set in $\mathbb{R}^2$ and $x\in \Omega$. Then for  any harmonic function $h$ on $\Omega$ and  $l<l_{x,\mathcal{P}}$, we have
\begin{equation}
|I(h,x,l)-\sum_{p\in\mathcal{P}}A_pI_p(h,x,l)|\leq \sum_{k=1}^{\infty}\frac{1}{(4k+2)!}\bigg\arrowvert\frac{\partial^{4k}h}{\partial\xi_1^{4k}}(x)\bigg\arrowvert\Big(2^{2k+1}+\|\mathcal{P}\|^{4k+2}\Big)(\frac{l}{2})^{4k},
\end{equation}
where
\begin{equation}
\|\mathcal{P}\|=\max_{p\in\mathcal{P}}\|2p+(1,1)\|,
\end{equation}
and $l_{x,\mathcal{P}}=min_{p\in\mathcal{P}}l_{x,p}$ with $l_{x,p}$ being the same as in Lemma 2.3.}

\noindent\textit{Proof.} Noticing  that
\begin{equation*}
I(h,x,l)=I_\theta(h,x,l), \text{ with } \theta=(0,0),
\end{equation*}
by using (2.10), we obtain that
\begin{equation*}
\begin{aligned}
&|I(h,x,l)-\sum_{p\in\mathcal{P}}A_pI_p(h,x,l)|\\
= \ &\bigg\arrowvert\big(1-8\sum_{p\in\mathcal{P}}c_pA_p\big)h(x)+\sum_{k=1}^{\infty}\frac{1}{(4k+2)!}\frac{\partial^{4k}h}{\partial\xi_1^{4k}}(x)\Big(T_\theta^{(k)}-\sum_{p\in\mathcal{P}}A_pT_p^{(k)}\Big)(\frac{l}{2})^{4k}\bigg\arrowvert\\
= \ 
&\sum_{k=1}^{\infty}\frac{1}{(4k+2)!}\bigg\arrowvert\frac{\partial^{4k}h}{\partial\xi_1^{4k}}(x)\bigg\arrowvert\Big\arrowvert T_\theta^{(k)}-\sum_{p\in\mathcal{P}}A_pT_p^{(k)}\Big\arrowvert(\frac{l}{2})^{4k}.
\end{aligned}
\end{equation*}

By calculation $T_\theta^{(k)}=(-1)^k2^{2k+1}$ and using the estimate (2.9) we then have
\begin{equation*}
\begin{aligned}
\Big\arrowvert &T_\theta^{(k)}-\sum_{p\in\mathcal{P}}A_pT_p^{(k)}\Big\arrowvert\\
\leq \  &2^{2k+1}+\sum_{p\in\mathcal{P}}8c_pA_p\|2p+(1,1)\|^{4k+2}\\
\leq \  &2^{2k+1}+\|\mathcal{P}\|^{4k+2},
\end{aligned}
\end{equation*}
which gives (2.15).
\hfill$\square$

\noindent\textbf{Remark.} In fact, given $\mathcal{P}$, by choosing $A_p$ properly, we can get ``higher'' rate of convergence for (2.15). Define
\begin{equation}
N=\inf\{k\in\mathbb{N}|\sum_{p\in\mathcal{P}}A_pT_p^{(k)}\ne T_\theta^{(k)}\},
\end{equation}
called the \textit{mean value level}  of $\mathcal{P}$ with coefficients $\{A_p\}_{p\in\mathcal{P}}$. Then (2.15) becomes
\begin{equation}
|I(h,x,l)-\sum_{p\in\mathcal{P}}A_pI_p(h,x,l)|\leq \sum_{k=N}^{\infty}\frac{1}{(4k+2)!}\bigg\arrowvert\frac{\partial^{4k}h}{\partial\xi_1^{4k}}(x)\bigg\arrowvert\Big(2^{2k+1}+\|\mathcal{P}\|^{4k+2}\Big)(\frac{l}{2})^{4k}.
\end{equation}
We will discuss more on the mean value level in Section 4.

\section{The Resistance form and the Laplacian on the unit square $S$}

In this section, we will show the expressions of the \textit{resistance form} and the \textit{(symmetric) Laplacian} on the unit square $S$ in terms of average values on cells.

\noindent\textbf{Lemma 3.1.} \textit{Let $\mathcal{P}$ and $\{A_p\}_ {p\in\mathcal{P}}$ be defined as in Theorem 2.4,  $\Omega$ be an open set in $\mathbb{R}^2$. Then for any $f\in C^1(\Omega)$ and $x\in\Omega$, we have
\begin{equation}
|\nabla f(x)|^2=\frac{1}{2}\mathcal{M}_\mathcal{P}\lim_{l\to 0}\frac{1}{l^2} \sum_{p\in\mathcal{P}}\bigg(A_p\!\!\sum_{\tiny\substack{
\mathcal{D}(x',l) \\ 
\sim_p \mathcal{D}(x,l)}
}\!\!\big(I(f,x',l)-I(f,x,l)\big)^2\bigg),
\end{equation}
where
\begin{equation}
\mathcal{M}_\mathcal{P}=\frac{1}{2}\Big(\sum_{p\in\mathcal{P}}\|p\|^2c_pA_p\Big)^{-1}.
\end{equation}}
\noindent\textit{Proof.} Let $l<l_{x,\mathcal{P}}$ where $l_{x,\mathcal{P}}$ is the same as defined in Theorem 2.4. Then for each $p\in\mathcal{P}$, there are $8c_p$ $p$-neighbors of the $l$-square $\mathcal{D}(x,l)$. It is easy to check that $\mathcal{D}(x',l)$ with $x'=(x_1+p_1l, x_2+p_2l)$ is one of them. Using the mean value theorem for intergral we have
\begin{equation*}
\begin{aligned}
I(f,x',l)-I(f,x,l)&=\frac{1}{l^2}\int_{\mathcal{D}(x',l)}f(\xi_1,\xi_2)d\mu(\xi)-\frac{1}{l^2}\int_{\mathcal{D}(x,l)}f(\xi_1,\xi_2)d\mu(\xi)\\
&=\frac{1}{l^2}\int_{\mathcal{D}(x,l)}\Big(f(\xi_1+p_1l,\xi_2+p_2l)-f(\xi_1,\xi_2)\Big)d\mu(\xi)\\
&=f(\eta_1+p_1l,\eta_2+p_2l)-f(\eta_1,\eta_2)
\end{aligned}
\end{equation*}
for some $(\eta_1,\eta_2)\in\mathcal{D}(x,l)$. Hence
\begin{equation*}
\lim_{l\to 0}\frac{1}{l}\big(I(f,x',l)-I(f,x,l)\big)=p_1\frac{\partial f}{\partial\xi_1}(x)+p_2\frac{\partial f}{\partial\xi_2}(x).
\end{equation*}

Dealing with other $p$-neighbors similarly, and summing over all the $8c_p$ terms, we get
\begin{equation*}
\begin{aligned}
&\lim_{l\to 0}\!\!\sum_{\tiny\substack{
\mathcal{D}(x',l) \\ 
\sim_p \mathcal{D}(x,l)}
}\!\!\frac{1}{l^2}\big(I(f,x',l)-I(f,x,l)\big)^2\\
=& \ c_p\Bigg(\Big(p_1\frac{\partial f}{\partial\xi_1}(x)+p_2\frac{\partial f}{\partial\xi_2}(x)\Big)^2+\Big(p_2\frac{\partial f}{\partial\xi_1}(x)+p_1\frac{\partial f}{\partial\xi_2}(x)\Big)^2+\Big(-p_1\frac{\partial f}{\partial\xi_1}(x)+p_2\frac{\partial f}{\partial\xi_2}(x)\Big)^2\\
&+\Big(-p_2\frac{\partial f}{\partial\xi_1}(x)+p_1\frac{\partial f}{\partial\xi_2}(x)\Big)^2+\Big(-p_1\frac{\partial f}{\partial\xi_1}(x)-p_2\frac{\partial f}{\partial\xi_2}(x)\Big)^2+\Big(-p_2\frac{\partial f}{\partial\xi_1}(x)-p_1\frac{\partial f}{\partial\xi_2}(x)\Big)^2 
\end{aligned}
\end{equation*}
\begin{equation*}
\begin{aligned}
&+\Big(p_1\frac{\partial f}{\partial\xi_1}(x)-p_2\frac{\partial f}{\partial\xi_2}(x)\Big)^2+\Big(p_2\frac{\partial f}{\partial\xi_1}(x)-p_1\frac{\partial f}{\partial\xi_2}(x)\Big)^2\Bigg) \qquad\qquad\qquad\qquad\qquad\qquad\qquad\\
=& \ 4\|p\|^2c_p|\nabla f(x)|^2.
\end{aligned}
\end{equation*}

Summing over all $p\in\mathcal{P}$ with the coefficient $A_p$, we obtain (3.1).\hfill$\square$

Now we turn to the resistance form on the unit square $S$. Notice that the Lebesgue measure $\mu$ restricted to $S$ becomes a regular probability measure with equal weights. Analogous to the $SG$ case, for a finite word $w=(w_1,\cdots,w_m)$ with each $w_j\in\{1,2,3,4\}$, define the average value for a function $f$ on $F_wS$ as
$$B_w(f)=\frac{1}{\mu(F_wS)}\int_{F_wS}fd\mu.$$

\noindent\textbf{Theorem 3.2.}  \textit{Let $\mathcal{P}$ and $\{A_p\}_ {p\in\mathcal{P}}$ be defined as before. For any $f,g\in C^1(S)$, $m\geq 0$, define
\begin{equation}
\mathcal{E}_m(f,g)=\mathcal{M}_\mathcal{P} \sum_{p\in\mathcal{P}}\bigg(A_p\!\!\sum_{\tiny\substack{
|w|=|w'|=m,  \\ 
F_{w'}S\sim_pF_wS}
}\!\!\big(B_{w'}(f)-B_w(f)\big)\big(B_{w'}(g)-B_w(g)\big)\bigg),
\end{equation}
where $\mathcal{M}_\mathcal{P}$ is the same as (3.2). Then we have
\begin{equation}
\mathcal{E}(f,g):=\lim_{m\to \infty}\mathcal{E}_m(f,g)=\int_S\nabla f\cdot\nabla gd\mu.
\end{equation}}

\noindent\textit{Proof.}
For a $m$-cell $F_wS$ in the unit square $S$, let $l_m$ denote the side length and $x_w$ the center of $F_wS$. Then 
\begin{equation*}
l_m=\frac{1}{2^m}, \ \mu(F_wS)=\frac{1}{4^m}, \text{ and } x_w=F_w\Big(\frac{1}{4}\sum_{i}q_i\Big),
\end{equation*}
hence
\begin{equation*}
B_w(f)=I(f,x_w,l_m), \ F_wS=\mathcal{D}(x_w,l_m) .
\end{equation*}

 Noticing that for $p\in \mathcal{P}$, the $p$-neighbors of $F_wS$ may not be within $S$, we define
\begin{equation*}
S_m=\bigcup\{F_wS:   |w|=m \text{  and all } p \text{-neighbors of } F_wS \text{ are contained in } S\},
\end{equation*}
and $S_m^c=S\setminus S_m$.
Then we have
\begin{equation*}
\begin{aligned}
\mathcal{E}_m(f,f)= \ &\mathcal{M}_\mathcal{P} \sum_{p\in\mathcal{P}}\bigg(A_p\!\!\sum_{\tiny\substack{
|w|=|w'|=m,  \\ 
F_{w'}S\sim_pF_wS}
}\!\!\big(B_{w'}(f)-B_w(f)\big)^2\bigg)\\
= \ &\frac{1}{2}\mathcal{M}_\mathcal{P}\sum_{|w|=m} \sum_{p\in\mathcal{P}}\bigg(A_p\!\!\sum_{\tiny\substack{
F_{w'}S \\ 
\sim_pF_wS}
}\!\!\big(B_{w'}(f)-B_w(f)\big)^2\bigg)\\
= \ &\frac{1}{2}\mathcal{M}_\mathcal{P}\!\!\sum_{\substack{
|w|=m \\
F_wS\subseteq S_m}
}\!\!\Bigg(\frac{1}{l_m^2}\sum_{p\in\mathcal{P}}\bigg(A_p\!\!\sum_{\tiny\substack{
\mathcal{D}(x_{w'},l_m) \\ 
\sim_p \mathcal{D}(x_w,l_m)}
}\!\!\big(I(f,x_{w'},l_m)-I(f,x_w,l_m)\big)^2\bigg)\mu(F_wS)\Bigg)
\end{aligned}
\end{equation*}
\begin{equation*}
\quad\quad\quad\quad\quad+\frac{1}{2}\mathcal{M}_\mathcal{P}\!\!\sum_{\substack{
|w|=m \\
F_wS\subseteq S_m^c}
}\!\!\Bigg(\frac{1}{l_m^2}\sum_{p\in\mathcal{P}}\bigg(A_p\!\!\sum_{\tiny\substack{
\mathcal{D}(x_{w'},l_m) \\ 
\sim_p \mathcal{D}(x_w,l_m)}
}\!\!\big(I(f,x_{w'},l_m)-I(f,x_w,l_m)\big)^2\bigg)\mu(F_wS)\Bigg).
\end{equation*}
It is obvious that $\lim_{m\to\infty}\mu(S_m^c)=0$. So by Lemma 3.1 we know that when  $m$ goes to $\infty$, the first term converges to $\int_S|\nabla f|^2d\mu$ and the second term converges to $0$. Hence
\begin{equation*}
\mathcal{E}(f,f)=\lim_{m\to\infty}\mathcal{E}_m(f,f)=\int_S|\nabla f|^2d\mu.
\end{equation*}
Then by using the polarization identity
\begin{equation}
\mathcal{E}_m(f,g)=\frac{1}{4}\Big(\mathcal{E}_m(f+g,f+g)-\mathcal{E}_m(f-g,f-g)\Big),
\end{equation}
we obtain
\begin{equation*}
\begin{aligned}
\mathcal{E}(f,g)=&\lim_{m\to\infty}\mathcal{E}_m(f,g)\\
=& \ \frac{1}{4}\bigg(\int_K|\nabla (f+g)|^2d\mu-\int_K|\nabla (f-g)|^2d\mu\bigg)\\
=& \ \int_K\nabla f\cdot\nabla gd\mu.
\end{aligned}
\end{equation*}
\hfill$\square$

We should remark that (3.3) and (3.4) imply that the renormalization factor of the resistance form  equals to $1$ in this case.

Next we come to the Laplacian on $S$.

\noindent\textbf{Lemma 3.3.} \textit{Let $\mathcal{P}$ and $\{A_p\}_ {p\in\mathcal{P}}$ be defined as before. $\Omega$ be an open set in $\mathbb{R}^2$. Then for  any $f\in C^2(\Omega)$ and $x\in\Omega$, we have
\begin{equation}
\Delta f(x)=\mathcal{M}_\mathcal{P}\lim_{l\to 0}\frac{1}{l^2}\bigg(\sum_{p\in\mathcal{P}}A_pI_p(f,x,l)-I(f,x,l)\bigg),
\end{equation}
where $\mathcal{M}_\mathcal{P}$ is the same as (3.2).}

\noindent\textit{Proof.} Since $8\sum_{p\in\mathcal{P}}c_pA_p=1$, we have
\begin{equation*}
\begin{aligned}
\sum_{p\in\mathcal{P}}A_pI_p(f,x,l)-I(f,x,l)=& \ \sum_{p\in\mathcal{P}}A_p\Big(I_p(f,x,l)-8c_pI(f,x,l)\Big)\\
=& \ \sum_{p\in\mathcal{P}}A_p\bigg(\!\!\sum_{\tiny\substack{
\mathcal{D}(x',l) \\ 
\sim_p \mathcal{D}(x,l)}
}\!\!\big(I(f,x',l)-I(f,x,l)\big)\bigg).
\end{aligned}
\end{equation*}
Let $p\in\mathcal{P}$. As in the proof of Lemma 3.1, there are $8c_p$ $p$-neighbors  of the $l$-square $\mathcal{D}(x,l)$, for example, $\mathcal{D}(x',l)$ and $\mathcal{D}(x'',l)$  with $x'=(x_1+p_1l, x_2+p_2l)$ and  $x''=(x_1-p_1l, x_2-p_2l)$ are two of them. By the mean value theorem for integral, we have
\begin{equation*}
\begin{aligned}
& \ I\big(f,x',l\big)-I\big(f,x,l\big)+I\big(f,x'',l\big)-I\big(f,x,l\big)\\
=& \ f(\eta_1+p_1l,\eta_2+p_2l)-f(\eta_1,\eta_2)+f(\eta_1-p_1l,\eta_2-p_2l)-f(\eta_1,\eta_2), \text{ for some }  (\eta_1,\eta_2)\in\mathcal{D}(x,l).
\end{aligned}
\end{equation*}
Thus
\begin{equation*}
\begin{aligned}
\lim_{l\to 0}\frac{1}{l^2}\Big(I\big(f,x',l\big)-I\big(f,x,l\big)+I\big(f,x'',l\big)-I\big(f,x,l\big)\Big)=p_1^2\frac{\partial^2f}{\partial\xi_1^2}(x)+2p_1p_2\frac{\partial^2f}{\partial\xi_1\partial\xi_2}(x)+p_2^2\frac{\partial^2f}{\partial\xi_2^2}(x).
\end{aligned}
\end{equation*}
Dealing with other $p$-neighbors similarly, and summing over all of them, we obtain
\begin{equation*}
\begin{aligned}
&\lim_{l\to 0}\frac{1}{l^2}\!\!\sum_{\tiny\substack{
\mathcal{D}(x',l) \\ 
\sim_p \mathcal{D}(x,l)}
}\!\!\!\!\big(I(f,x',l)-I(f,x,l)\big)\\
=& \ c_p\bigg(p_1^2\frac{\partial^2f}{\partial\xi_1^2}(x)+2p_1p_2\frac{\partial^2f}{\partial\xi_1\partial\xi_2}(x)+p_2^2\frac{\partial^2f}{\partial\xi_2^2}(x)+p_2^2\frac{\partial^2f}{\partial\xi_1^2}(x)+2p_1p_2\frac{\partial^2f}{\partial\xi_1\partial\xi_2}(x)+p_1^2\frac{\partial^2f}{\partial\xi_2^2}(x)\\
&+(-p_1)^2\frac{\partial^2f}{\partial\xi_1^2}(x)+2(-p_1)p_2\frac{\partial^2f}{\partial\xi_1\partial\xi_2}(x)+p_2^2\frac{\partial^2f}{\partial\xi_2^2}(x)+(-p_2)^2\frac{\partial^2f}{\partial\xi_1^2}(x)+2p_1(-p_2)\frac{\partial^2f}{\partial\xi_1\partial\xi_2}(x)\\
&+p_1^2\frac{\partial^2f}{\partial\xi_2^2}(x)\bigg)\\
=& \ 2\|p\|^2c_p\Delta f(x).
\end{aligned}
\end{equation*}
Summing over all $p\in\mathcal{P}$, we obtain (3.6).\hfill$\square$

Then we could apply Lemma 3.3 to the unit square $S$ to reconstruct the Laplacian in a self-similar manner.

\noindent\textbf{Theorem 3.4.} \textit{Let $\mathcal{P}$ and $\{A_p\}_ {p\in\mathcal{P}}$ be defined as before. For
any $f\in C^2(S)$, $x\in S\setminus\partial S$, $m\geq 0$, define
\begin{equation}
\Delta_mf(x)=\mathcal{M}_\mathcal{P}4^m\bigg(\sum_{p\in\mathcal{P}}\big(A_p\!\!\sum_{\tiny\substack{
F_{w'}S \\ 
\sim_pF_wS}
}\!\!B_{w'}(f)\big)-B_w(f)\bigg) \text{ for } x\in F_wS\text{ with } |w|=m,
\end{equation}
where $\mathcal{M}_\mathcal{P}$ is the same as (3.2). Then
\begin{equation}
\Delta f(x)=\lim_{m\to\infty}\Delta_mf(x)
\end{equation} uniformly.}

\noindent\textit{Proof.} Similarly, we define $S_m$ and $S_m^c$ as we did in the proof of Theorem 3.2. Since $x\in S\setminus\partial S$ we know that there exists an integer $m_0$ such that when $m\geq m_0$ (that is, $l_m\leq 1/2^{m_0}$),  $x\in F_wS\subseteq S_m$ for some $w$ of length $m$. Obviously, $x_w$(the center of $F_wS$) will goes to $x$ as $m$ goes to $\infty$. Then from (3.6) we get
\begin{equation*}
\begin{aligned}
\lim_{m\to\infty}\Delta_mf(x)=& \ \mathcal{M}_\mathcal{P}\lim_{m\to\infty}4^m\bigg(\sum_{p\in\mathcal{P}}\big(A_p\!\!\sum_{\tiny\substack{
F_{w'}S \\ 
\sim_pF_wS}
}\!\!B_{w'}(f)\big)-B_w(f)\bigg)\\
=& \ \mathcal{M}_\mathcal{P}\lim_{l_m\to 0}\frac{1}{l_m^2}\bigg(\sum_{p\in\mathcal{P}}A_pI_p(f,x_w,l_m)-I(f,x_w,l_m)\bigg)\\
=& \ \Delta f(x).
\end{aligned}
\end{equation*}
The uniform convergence comes from the fact that $S$ is compact.
\hfill$\square$

\section{The level of mean value property}

In Section 2 we have seen that the rate of convergence of (2.15) is decided by the mean value level of $(\mathcal{P}, \{A_p\}_{p\in \mathcal{P}})$. For given $\mathcal{P}$ and  positive integer $N$, from (2.17), we may find $\{A_p\}_{p\in\mathcal{P}}$ by solving equations

\begin{equation}
\begin{aligned}
8\sum_{p\in\mathcal{P}}c_pA_p&=1,\\
\sum_{p\in\mathcal{P}}A_pT_p^{(k)}&=T_\theta^{(k)}, \text{ for } k=0,1,\dots,N\!-\!1,
\end{aligned}
\end{equation}
such that $N$ is the mean value level of $(\mathcal{P},\{A_p\}_{p\in\mathcal{P}})$.

Here are some solutions:

For $\mathcal{P}=\{(0,1)\},  N=1$, we have a unique solution
\begin{equation*}
A_{(0,1)}=\frac{1}{4}, \
\mathcal{M}_\mathcal{P}=4;
\end{equation*}

For $\mathcal{P}=\{(0,1),(1,1)\},  N=1$, we have infinite solutions satisfying
\begin{equation*}
A_{(0,1)}+A_{(1,1)}=\frac{1}{4};
\end{equation*}

For $\mathcal{P}=\{(0,1),(1,1)\},  N=2$, we have a unique solution
\begin{equation*}
A_{(0,1)}=\frac{1}{5}, \ A_{(1,1)}=\frac{1}{20}, \ 
\mathcal{M}_\mathcal{P}=\frac{10}{3};
\end{equation*}

For $\mathcal{P}=\{(0,1),(1,1),(0,2)\},  N=3$, we have a unique solution
\begin{equation*}
A_{(0,1)}=\frac{16}{75} , \ A_{(1,1)}=\frac{1}{25}, \ A_{(0,2)}=-\frac{1}{300}, \ 
\mathcal{M}_\mathcal{P}=\frac{25}{7};
\end{equation*}

For $\mathcal{P}=\{(0,1),(1,1),(0,2),(1,2)\}, N=4$, we have a unique solution
\begin{equation*}
A_{(0,1)}=\frac{38}{183}, \ A_{(1,1)}=\frac{103}{2379}, \ A_{(0,2)}=-\frac{17}{9516}, \ A_{(1,2)}=\frac{1}{2379}, \ 
\mathcal{M}_\mathcal{P}=\frac{793}{231}.
\end{equation*}

On the other hand, for a given $\mathcal{P}$ with coefficients $\{A_p\}_{p\in\mathcal{P}}$, it is natural to ask which kind of  harmonic functions does satisfy the mean value property exactly.

\noindent\textbf{Theorem 4.1} \textit{Let $(\mathcal{P},\{A_p\}_{p\in\mathcal{P}})$ be defined as before with the mean value level $N$, $\Omega$ be a connected open set in $\mathbb{R}^2$, and  $h$ be a harmonic function on $\Omega$. Suppose there is a connected open subset $U$ of $\Omega$ such that  for any $x\in U$, the identity 
\begin{equation*}
I(h,x,l)-\sum_{p\in\mathcal{P}}A_pI_p(h,x,l)=0
\end{equation*}
holds for sufficiently small $l$. Then $h$ is a polynomial harmonic function on $\Omega$ with degree no more than $4N$.}

\noindent\textit{Proof.} By applying (2.10) and (2.17), we know that for any $x\in U$ and sufficiently small $l$,
\begin{equation*}
I(h,x,l)-\sum_{p\in\mathcal{P}}A_pI_p(h,x,l)=\sum_{k=N}^{\infty}\frac{1}{(4k+2)!}\frac{\partial^{4k}h}{\partial\xi_1^{4k}}(x)\Big(T_\theta^{(k)}-\sum_{p\in\mathcal{P}}A_pT_p^{(k)}\Big)(\frac{l}{2})^{4k}=0.
\end{equation*}
From the arbitrariness of $l$, we know that for any $k\geq N$,
\begin{equation*}
\frac{\partial^{4k}h}{\partial\xi_1^{4k}}(x)\Big(T_\theta^{(k)}-\sum_{p\in\mathcal{P}}A_pT_p^{(k)}\Big)=0.
\end{equation*}
Then by the definition of $N$,
\begin{equation*}
T_\theta^{(N)}-\sum_{p\in\mathcal{P}}A_pT_p^{(N)}\ne 0,
\end{equation*}
thus
\begin{equation*}
\frac{\partial^{4N}h}{\partial\xi_1^{4N}}(x)=0,  \forall x\in U.
\end{equation*}
Hence from (2.12), we can expand $h$ near $x$ as a polynomial harmonic function of degree no more than $4N$ as follows
\begin{equation*}
\begin{aligned}
h(\xi)=&h(x)+\sum_{n=1}^{4N-1}\bigg(\frac{\partial^nh}{\partial\xi_1^n}(x)\sum_{j=0}^{[\frac{n}{2}]}(-1)^j\frac{(\xi_1-x_1)^{n-2j}}{(n-2j)!}\frac{(\xi_2-x_2)^{2j}}{(2j)!}\bigg)\\
&+\sum_{n=1}^{4N}\bigg(\frac{\partial^nh}{\partial\xi_1^{n-1}\partial\xi_2}(x)\sum_{j=0}^{[\frac{n-1}{2}]}(-1)^j\frac{(\xi_1-x_1)^{n-2j-1}}{(n-2j-1)!}\frac{(\xi_2-x_2)^{2j+1}}{(2j+1)!}\bigg).
\end{aligned}
\end{equation*}
Since $h$ is harmonic on $\Omega$, $h$ is real analytic in $\Omega$. Then from the principle of analytic continuation, we know that the above identity is valid for any $x\in\Omega$. \hfill$\square$

\end{document}